\documentclass[12pt]{article} 

\usepackage{amsmath,amsxtra,amssymb,latexsym, amscd}
\usepackage[mathscr]{eucal}

\begin{document}  
\title{\bf Spectrum of a bounded sequence and inhomogeneous delay linear difference equations  in
a Banach space}

\def\1{\rule{0cm}{0cm}} \def\qd{\rule{3mm}{3mm}} \def\BB{$\bullet$}
\renewcommand{\arraystretch}{1.25}
\renewcommand{\theequation}{\thesectn.\arabic{equation}}
\def\sce{\setcounter{equation}{0}}  \newcounter{sectn} 
\newcounter{sbsect}
\def\sect#1{\addtocounter{section}{1}\sce\setcounter{sbsect}{0}%
        \renewcommand{\thesectn}{\thesection}\1\smallskip\\
        {\1\hspace{-2em}\large\bf\thesectn.\qquad #1\smallskip\par}}
\def\subsect#1{\addtocounter{sbsect}{1}\sce%
        
\renewcommand{\thesectn}{\thesection:\Alph{sbsect}}\1\smallskip\\
        {\bf\1\hspace{-1.5em}\thesectn.\qquad #1\smallskip\par}}
\newtheorem{Theorem}{THEOREM} \newtheorem{Lemma}[Theorem]{LEMMA}
\newtheorem{Corollary}[Theorem]{COROLLARY}
\def\thm#1#2{\be{Theorem}{\lb{#1} #2}} \def\LEM#1#2{\BE{Lemma}{\LB{#1} 
#2}}
\def\COR#1#2{\BE{Corollary}{\LB{#1} #2}}
\def\proof{\bigskip\noindent {\sc Proof:}\qquad}
\def\REM{\1\smallskip\par\noindent{\bf REMARK:}\qquad }
\def\qed{\hfill$\quad$\qd\medskip\\} \def\ds{\displaystyle}
\def\LB#1{\label{#1}} \def\BE#1#2{\begin{#1} #2 \end{#1}}
\def\EQ#1#2{\BE{equation}{\LB{#1} #2}} \def\ARR#1#2{\BE{array}{{#1} 
#2}}
\def\DES#1{\BE{description}{#1}} \def\QT#1{\BE{quote}{#1}}
\def\ENUM#1{\BE{enumerate}{#1}} \def\ITM#1{\BE{itemize}{#1}}
 \def\COM#1{\par\noindent{\bf COMMENT:\quad\sl #1}\par\noindent}
\def\mapsfrom{\hbox{$\;{\leftarrow}\kern-.15em{\mapstochar}\:\:$}}
\def\vv{\kern.344em{\rule[.18ex]{.075em}{1.32ex}}\kern-.344em}
\def\RE{\mbox{\rm I\kern-.21em R}} \def\CX{\mbox{\rm \vv C}}
\def\imp{\Rightarrow} \def\emb{\hookrightarrow} 
\def\wk{\rightharpoonup}
\def\rd{\dot{\1}} \def\d{\cdot} \def\+{\oplus} \def\x{\times}
\def\<{\langle} \def\>{\rangle} \def\o{\circ} \def\at#1{\Bigr|_{#1}}
\def\cd{\partial} \def\grad{\nabla} \def\L{\left} \def\R{\right}
\def\A{\mathbf  A} \def\by{\mathbf{y}} \def\bS{\mathbf{S}}
\def\I{{\cal I}} \def\A{{\mathbf A}} \def\D{{\mathbb D}}\def\bc{\mathbf{c}}
\def\X{\mathbf{X}}
\def\Y{\mathbb{Y}}
 \def\bo{\mathbf{O}}
 \def\B{\mathbf{B}}\def\bone{\mathbf{1}}\def\bzero{\mathbf{0}}
\def\H{{\mathcal H}} \def\U{{\cal U}} \def\bc{\mathbf{c}}
\def\F{\mathbf{F}}
\def\eq{equation} \def\de{differential \eq} \def\pde{partial \de}
\def\sol{solution} \def\pb{problem} \def\bdy{boundary} 
\def\fn{function}
\def\dde{delay \de} \def\ev{eigenvalue}
\def\R{\mathbb R}
\def\Q{\mathbb Q}
\def\C{\mathbb C}
\def\N{\mathbb N}
\def\Z{\mathbb Z}
\def\X{\mathbb X}
\author{
{\bf Dang Vu Giang}\\
Hanoi Institute of Mathematics\\
Vietnam Academy of Science and Technology\\
18 Hoang Quoc Viet, 10307 Hanoi, Vietnam\\
{\footnotesize          e-mail: $\<$dangvugiang@yahoo.com$\>$}\\
\1\\
}
\maketitle
                      
\noindent {\bf Abstract.}    We  study the asymptotic behavior  of a bounded solution of an inhomogeneous delay linear difference equation  in a Banach space by using the spectrum of  bounded sequences. 
We get a significant extension of excellent results in [1]. A new simple proof is also found for the famous Gelfand spectral radius theorem. Moreover, among other things we  prove that if the spectrum of a bounded sequence $\{x_n\}_n$ is finite then $x_n=c_1\vartheta_1^n+c_2\vartheta_2^n+\cdots+c_k\vartheta_k^n+o(1)$ as $n\to\infty$ where  $|\vartheta_1|=|\vartheta_2|=\cdots=|\vartheta_k|=1$.
 


\eject

\sect{Introduction} 
\bigskip
\par\noindent 
Katznelson and Tzafriri \cite{Katznelson} proved the following famous result.

\bigskip
\noindent {\bf Theorem A. } {\it  Let $T:\mathbb{X}\to \mathbb{X}$ denote  a bounded linear operator and $\partial \mathbb{D}$ the unit circle. If $T$ is power bounded (that is the sequence of norms $\left\{ \left\| T \right\|,\left\| {{T}^{2}} \right\|,\cdots  \right\}$ is bounded) and $\partial \mathbb{D}\cap \sigma \left( T \right)\subseteq \left\{ 1  \right\}$ then
 \[\underset{n\to \infty }{\mathop{\lim }}\,\left( {{T}^{n+1}}- {{T}^{n}} \right)=0.\]}

\medskip
\noindent Vu Quoc Phong \cite{Vu} reproved this theorem. In this paper, we use a new method (the spectrum of a bounded sequence in a Banach space)  to prove that there exists the limit of $T^n$ as $n\to\infty$ ($T$ denotes the operator in Theorem A). 

\sect{ 	Holomorphic functions in a Banach space and Resolvent of an isometry operator}

\bigskip
\par\noindent 
Let $\mathbb{X}=\left( \mathbb{X},{{\left\| {\cdot} \right\|}_{\mathbb{X}}} \right)$ denote a Banach space. We are interested in those functions $f:\mathbb{C}\to \mathbb{X}$ which can be expressed as power series $f\left( z \right)=\sum\limits_{k=0}^{\infty }{{{z}^{k}}{{x}_{k}}}$ where $\left\{ {{x}_{k}} \right\}_{k=0}^{\infty }\subseteq \mathbb{X}.$ This series is convergent in the norm topology of $\mathbb{X}.$ This means that 
$\sum\limits_{k=0}^{\infty }{{{\left\| {{x}_{k}} \right\|}_{\mathbb{X}}}{{\left| z \right|}^{k}}}<\infty $ for every $z\in \mathbb{C}.$ These functions are called entire functions in the Banach space $\mathbb{X}.$ If this series is finite we say about polynomial function (operator).
The complex derivative of $f$ is ${f}'\left( z \right)=\sum\limits_{k=1}^{\infty }{k{{z}^{k-1}}{{x}_{k}}.}$  Moreover, the complex integral of $f$ gives
\[{{x}_{k}}=\frac{1}{2\pi i}\oint\limits_{\left| z \right|=R}{\frac{f\left( z \right)dz}{{{z}^{k+1}}}}\quad \text{  for } k=0,1,2,\cdots .\]
Hence, if  ${{\left\| f\left( z \right) \right\|}_{\mathbb{X}}}$  is bounded then $f$ is constant (a vector of the Banach space $\mathbb{X}$).
We remind the readers that the path integration of a continuous function $\varphi :\mathbb{C}\to \mathbb{X}$ along the simple piecewise differentiable curve $\gamma :\left[ a,b \right]\to \mathbb{C}$ is 
\[\int\limits_{\gamma }{\varphi \left( z \right)dz}=\int\limits_{a}^{b}{{\gamma }'\left( t \right)\varphi \left( \gamma \left( t \right) \right)dt.}\]
Now let $\Omega $ be an open (unbounded)  region of complex plane. Consider those functions $g:\Omega \to \mathbb{X}$ such that for any ${{z}_{0}}\in \Omega $ there is $\delta >0$ such that
$g\left( z \right)=\sum\limits_{k=0}^{\infty }{{{\left( z-{{z}_{0}} \right)}^{k}}{{x}_{k}}}$ for  some $\left\{ {{x}_{k}} \right\}_{k=0}^{\infty }\subseteq \mathbb{X}$ and $\left| z-{{z}_{0}} \right|<\delta .$
These functions of this condition are called holomorphic functions in the region $\Omega $. 
If we can extend $g$ to the whole complex plane without breaking this condition then $g$ is also called an entire function. Otherwise, we say about the essential singularity of $g.$ A point ${{z}_{0}}\in \mathbb{C}$ is called an essential singularity point of $g$ if $g$ cannot be extended to $\Omega \cup \left\{ {{z}_{0}} \right\}$ without breaking the holomorphy condition.
Most of time we are interested in the resolvent ${{\left( \lambda -A \right)}^{-1}}$ of a linear bounded operator $A:\mathbb{X}\to \mathbb{X}$. This is a holomorphic function defined on $\mathbb{C}\backslash \sigma \left( A \right)$ by Laurent series
\[{{\left( \lambda -A \right)}^{-1}}=\sum\limits_{n=0}^{\infty }{\frac{{{A}^{n}}}{{{\lambda }^{n+1}}}}\]
which is convergent for all $\left| \lambda  \right|>\rho \left( A \right)$  ($\rho(A)$ denotes the spectral radius of $A$).
On the other hand, if $\X$ is finite dimensional then the resolvent ${{\left( \lambda -A \right)}^{-1}}$ has finite poles in $\sigma(A)$ (the spectrum of $A$). Let $\chi_A$ denote the characteristic polynomial of $A$. Then $\chi_A(z)=0$ for every $z\in\sigma(A)$ and the resolvent ${{\left( \lambda -A \right)}^{-1}}$ is holomorphic in $\C\setminus\sigma(A)$. Therefore, we can write $${{\left( \lambda -A \right)}^{-1}}=\frac{\phi(\lambda)}{\chi_A(\lambda)},$$ where $\phi:\C\to B(\X)$ is a holomorphic function ($B(\X)$ denotes the set of continuous linear operators on $\X$ and $\chi_A$ denotes the characteristic polynomial of $A$). Now multiply with ${{\left( \lambda -A \right)}}\chi_A(\lambda)$ side by side we have $\chi_A(\lambda)I={{\left( \lambda -A \right)}}\phi(\lambda)$. Let $\lambda=A$ we have $\chi_A(A)=0$. The famous Caley-Hamilton theorem is proved. Moreover, we can prove the famous theorem of I. Gelfand on the spectral radius as follows. Let ${{a}_{n}}=\ln \left\| {{A}^{n}} \right\|.$ Then ${{a}_{n+m}}=\ln \left\| {{A}^{n+m}} \right\|\le \ln \left( \left\| {{A}^{n}} \right\|\left\| {{A}^{m}} \right\| \right)=\ln \left\| {{A}^{n}} \right\|+\ln \left\| {{A}^{m}} \right\|={{a}_{n}}+{{a}_{m}}$ and consequently, there is $\lim {{a}_{n}}/n=:\ln r,$ that is $\lim {{\left\| {{A}^{n}} \right\|}^{1/n}}=r.$ On the other hand, the relsovent series \[{{\left( \lambda -A \right)}^{-1}}=\sum\limits_{n=0}^{\infty }{\frac{{{A}^{n}}}{{{\lambda }^{n+1}}}}\]
 is absolutely convergent for all $\left| \lambda  \right|>\rho \left( A \right)$ 
and
\[\sum\limits_{n=0}^{\infty }{\frac{\left\| {{A}^{n}} \right\|}{\rho {{\left( A \right)}^{n}}}}=\infty .\]
 If  $\lim {{\left\| {{A}^{n}} \right\|}^{1/n}}=r<\rho \left( A \right)$ then $\left\| {{A}^{n}} \right\|<{{\left[ \rho \left( A \right)-\varepsilon  \right]}^{n}}$ for all $n>N$ and consequently, \[\infty =\sum\limits_{n>N}{\frac{\left\| {{A}^{n}} \right\|}{\rho {{\left( A \right)}^{n}}}}<\sum\limits_{n>N}{{{\left[ \frac{\rho \left( A \right)-\varepsilon }{\rho \left( A \right)} \right]}^{n}}}<\infty \]
 which is a contradiction. If $\lim {{\left\| {{A}^{n}} \right\|}^{1/n}}=r>\rho \left( A \right)$ then $\left\| {{A}^{n}} \right\|>{{\left[ \rho \left( A \right)+\varepsilon  \right]}^{n}}$ for all $n>N$ and    the series
\[\sum\limits_{n=0}^{\infty }{\frac{\left\| {{A}^{n}} \right\|}{{{\left| \lambda  \right|}^{n}}}}\]  divergent for $\left| \lambda  \right|=\rho \left( A \right)+\varepsilon >\rho \left( A \right)$ which is a contradiction. In the next sections
we are specially interested in the resolvent of isometry operators. More exactly, if $A:\mathbb{X}\to \mathbb{X}$ is an isometry linear operator then the spectrum of $A$ is contained in the unit circle and the norm of ${{\left( \lambda -A \right)}^{-1}}$ is bounded by 
${{\bigl| \left| \lambda  \right|-1 \bigr|}^{-1}}.$ On the other hand, any isolated essential singular point of ${{\left( \lambda -A \right)}^{-1}}$ is a simple pole \cite{Minh}.
\bigskip
\bigskip

\sect{ 	Spectrum of a bounded sequence}
\bigskip
\par\noindent 
Let $\mathbb{X}=\left( \mathbb{X},{{\left\| {\cdot} \right\|}_{\mathbb{X}}} \right)$ denote a Banach space. Let $\mathbf{x}=\left\{ {{x}_{0}},{{x}_{1}},\cdots  \right\}$ denote a sequence with elements in ${\mathbb X}$ and ${{\ell }^{\infty }}\left( \mathbb{X} \right)$ the Banach space of bounded sequences in ${\mathbb X}$ with the norm  $\left\| \mathbf{x} \right\|=\sup \left\{ {{\left\| {{x}_{0}} \right\|}_{\mathbb{X}}},{{\left\| {{x}_{1}} \right\|}_{\mathbb{X}}},\cdots  \right\}$. Moreover, let ${{c}_{0}}\left( \mathbb{X} \right)$  be the subpace of ${{\ell }^{\infty }}\left( \mathbb{X} \right)$ consisting of  vanishing sequences $\mathbf{x}=\left\{ {{x}_{0}},{{x}_{1}},\cdots  \right\}$ in ${\mathbb X}$ that is $\underset{n\to \infty }{\mathop{\lim }}\,{{x}_{n}}=0.$ Let  \[S:{{\ell }^{\infty }}\left( \mathbb{X} \right)\to {{\ell }^{\infty }}\left( \mathbb{X} \right)\] denote the shift operator, that is ${{\left( S\mathbf{x} \right)}_{n}}={{x}_{n+1}}$. Let \[\mathbb{Y}={{\ell }^{\infty }}\left( \mathbb{X} \right)/{{c}_{0}}\left( \mathbb{X} \right)\] be the quotient space. 
The equivalent class containing $\mathbf{x}=\left\{ {{x}_{0}},{{x}_{1}},\cdots  \right\}$ is denoted by $\mathbf{\bar{x}}=\left\{ {{{\bar{x}}}_{0}},{{{\bar{x}}}_{1}},\cdots  \right\}$
The norm of an $\mathbf{\bar x}=\overline{\left\{ {{x}_{0}},{{x}_{1}},\cdots  \right\}}\in\mathbb{Y}
$ is defined by $\left\|\mathbf{\bar x}\right\|_{\mathbb {Y}}
=\limsup_{n\to\infty}\left\|x_n
\right\|_{\X}$
 and the reduced  shift operator of  $S$  is denoted by $\overline{S}:{\mathbb{Y}}\to {\mathbb{Y}}$. Then $\overline{S}$ is an isometry operator so the spectrum of $\overline{S}$ is contained in the unit circle so the  resolvent operator $R\left( \lambda ,\overline{S} \right)$ of $\overline{S}$ is analytic and injective  for every $\left| \lambda  \right|\ne 1.$ 
Hence, if $R\left( \lambda ,\overline{S} \right)\mathbf{\bar{x}}=0$ for some $\left| \lambda  \right|\ne 1$ 
 then $\mathbf{\bar{x}}=0$ which means $\underset{n\to \infty }{\mathop{\lim }}\,{{x}_{n}}=0.$
 Moreover, the norm of $R\left( \lambda ,\overline{S} \right)$ is bounded by ${{\bigl| \left| \lambda  \right|-1 \bigr|}^{-1}}.$
 These conditions hold for resolvent of any isometry operator.
The spectrum of a bounded sequence $\mathbf{x}=\left\{ {{x}_{0}},{{x}_{1}},\cdots  \right\}$ denoted by $\sigma \left( \mathbf{x} \right)$ is the set of all  essential  (non-removable)   singular points of $g\left( \lambda  \right)=R\left( \lambda ,\overline{S} \right)\mathbf{\bar{x}}$
(holomorphic function taking values in $\Y$). Then $\sigma \left( \mathbf{x} \right)$ is contained in the unit circle $\partial \mathbb{D}$. Moreover, we have

\bigskip
\noindent {\bf Theorem 1. } {\it 
$\sigma \left( \mathbf{x} \right)$ is empty iff $\underset{n\to \infty }{\mathop{\lim }}\,{{x}_{n}}=0.$ }

\bigskip
\noindent {\sl Proof: }
Let $\mathbf{x}$ be a vanishing sequence. Then $\mathbf{\bar{x}}$ is zero and $g\left( \lambda  \right)=R\left( \lambda ,\overline{S} \right)\mathbf{\bar{x}}$ is identically 0. Therefore, $\sigma \left( \mathbf{x} \right)$ is empty. Now assume that $g\left( \lambda  \right)=R\left( \lambda ,\overline{S} \right)\mathbf{\bar{x}}$ is an entire function. 
Since the norm of $R\left( \lambda ,\overline{S} \right)$ is bounded by 
$\bigl| \left| \lambda  \right|-1 \bigr|^{-1}$,  
we have $\left\| g\left( \lambda  \right) \right\|_{\X}\le {{\bigl| \left| \lambda  \right|-1 \bigr|}^{-1}}\left\| {\bar{\mathbf x}} \right\|,$  
which means that $
\bigl|\left| \lambda  \right|-1 \bigr|\left\| g\left( \lambda  \right) \right\|_{\X}\le \left\| {\bar{\mathbf x}} \right\|$. 
Since $g\left( \lambda  \right) $ is an entire function, we should have 
$g\left( \lambda  \right) $ is identically 0 (bounded entire functions are constant). By the injectivity of the resolvent we should have 
$\mathbf{\bar{x}}=0$ which means $\underset{n\to \infty }{\mathop{\lim }}\,{{x}_{n}}=0.$ The proof is now complete.

\bigskip
\noindent {\bf Note 1. } Theorem 1 is presented in \cite{Minh} without strict proof. We refer \cite{Rudin} for readers intersted in complex function and spectral  theory .

\bigskip
\noindent {\bf Theorem 2. } {\it $\sigma \left( \mathbf{x} \right)=\left\{ 1 \right\}$ iff 
$\underset{n\to \infty }{\mathop{\lim }}\left({{x}_{n+1}}-{{x}_{n}}\right)=0.$
}

\bigskip
\noindent {\sl Proof: } If $\underset{n\to \infty }{\mathop{\lim }}\,\left( {{x}_{n+1}}-{{x}_{n}} \right)=0$  then $S\bar{\mathbf x}=\bar{\mathbf x}$ and $g\left( \lambda  \right)=R\left( \lambda ,\overline{S} \right)\mathbf{\bar{x}}=(\lambda -1)^{-1}{{\mathbf{\bar{\mathbf x}}}}$ so $\sigma \left( \mathbf{x} \right)=\left\{ 1 \right\}$. Now assume that $\sigma \left( \mathbf{x} \right)=\left\{ 1 \right\}$. Then 
$g\left( \lambda  \right)=R\left( \lambda ,\overline{S} \right)\mathbf{\bar{x}}={{\left( \lambda -1 \right)}^{-1}}\mathbf{\bar{c}}+\varphi \left( \lambda  \right)$ with $\varphi \left( \lambda  \right)$ is an entire function taking values in ${\mathbb X}$. Since $\left\| g\left( \lambda  \right) \right\|_{\X}\le {{\bigl| \left| \lambda  \right|-1 \bigr|}^{-1}}\left\| {\bar{\mathbf x}} \right\|,$ we deduce that 
$g\left( \lambda  \right)=\left( \lambda -1 \right)^{-1}\bar{\mathbf  c}$ and consequently,
$\lambda \left( \mathbf{\bar{c}}-\mathbf{\bar{x}} \right)=S\mathbf{\bar{c}}-\mathbf{\bar{x}}$ for every $\lambda \ne 1.$ Let $\lambda \to\infty $ we get  $\mathbf{\bar{c}}=\mathbf{\bar{x}}$ and
$S\bar{\mathbf x}=\bar{\mathbf x}$ 
so
$\underset{n\to \infty }{\mathop{\lim }}\left({{x}_{n+1}}-{{x}_{n}}\right)=0.$
 The proof is now complete.

\bigskip
\noindent {\bf Note 2. } Theorem 2 is also presented in \cite{Minh} without strict proof.  Similar to Theorem 2 we consider the case where $\sigma \left(\mathbf  x \right)$ has only one point. We have 

\bigskip
\noindent {\bf Theorem 3. } {\it
 $\sigma \left( \mathbf x \right)=\left\{ \vartheta  \right\}$  iff  $\underset{n\to \infty }{\mathop{\lim }}\,\left( {{x}_{n+1}}-\vartheta {{x}_{n}} \right)=0.$}

\bigskip
\noindent {\sl Proof: } If $\underset{n\to \infty }{\mathop{\lim }}\,\left( {{x}_{n+1}}-\vartheta{{x}_{n}} \right)=0$  then $S\bar{\mathbf x}=\vartheta\bar{\mathbf x}$ and $g\left( \lambda  \right)=R\left( \lambda ,\overline{S} \right)\mathbf{\bar{x}}=(\lambda -\vartheta)^{-1}{{\mathbf{\bar{\mathbf x}}}}$ so $\sigma \left( \mathbf{x} \right)=\left\{ \vartheta \right\}$. Now assume that $\sigma \left( \mathbf{x} \right)=\left\{ \vartheta \right\}$. Then 
$g\left( \lambda  \right)=R\left( \lambda ,\overline{S} \right)\mathbf{\bar{x}}={{\left( \lambda -\vartheta \right)}^{-1}}\mathbf{\bar{c}}+\varphi \left( \lambda  \right)
={{\left( \lambda\bar\vartheta-1 \right)}^{-1}\bar\vartheta}\mathbf{\bar{c}}+\varphi \left( \lambda  \right)
$ 
with $\varphi \left( \lambda  \right)$ is an entire function taking values in ${\mathbb X}$. Since $\left\| g\left( \lambda  \right) \right\|_{\X}\le {{\bigl| \left| \lambda  \right|-1 \bigr|}^{-1}}\left\| {\bar{\mathbf x}} \right\|,$ we deduce that 
$g\left( \lambda  \right)=\left( \lambda -\vartheta \right)^{-1}\bar{\mathbf  c}$ and consequently,
$\lambda \left( \mathbf{\bar{c}}-\mathbf{\bar{x}} \right)=S\mathbf{\bar{c}}-\vartheta\mathbf{\bar{x}}$ for every $\lambda \ne \vartheta.$ Let $\lambda \to\infty $ we get  $\mathbf{\bar{c}}=\mathbf{\bar{x}}$ and
$S\bar{\mathbf x}=\vartheta\bar{\mathbf x}$ 
so
$\underset{n\to \infty }{\mathop{\lim }}\left({{x}_{n+1}}-\vartheta{{x}_{n}}\right)=0.$
 The proof is now complete.

\bigskip
\noindent {\bf Note 3. } Theorems 1 and 3 give the following theorem which is invented by Katznelson and Tzafriri \cite{Katznelson} and reproved by Vu Quoc Phong \cite{Vu} in the case $\vartheta =1$.

\bigskip
\noindent {\bf Theorem 4. } {\it  Let $T:\mathbb{X}\to \mathbb{X}$ denote  a bounded linear operator and $\partial \mathbb{D}$ the unit circle. If $T$ is power bounded (that is the sequence of norms $\left\{ \left\| T \right\|,\left\| {{T}^{2}} \right\|,\cdots  \right\}$ is bounded) and $\partial \mathbb{D}\cap \sigma \left( T \right)\subseteq \left\{ \vartheta  \right\}$ then
 \[\underset{n\to \infty }{\mathop{\lim }}\,\left( {{T}^{n+1}}-\vartheta {{T}^{n}} \right)=0.\]
}

\bigskip
\noindent {\sl Proof: } Let ${{x}_{n}}={{T}^{n}}$ and consider the spectrum of  $\mathbf{x}=\left\{ {{x}_{0}},{{x}_{1}},\cdots  \right\}$. We have at once that $\sigma \left( \mathbf{x} \right)\subseteq \partial \mathbb{D}\cap \sigma \left( T \right)\subseteq \left\{ \vartheta  \right\}$ so by Theorems 1 and 3,  $\underset{n\to \infty }{\mathop{\lim }}\left({{x}_{n+1}}-\vartheta {{x}_{n}}\right)=0.$ The proof is now complete.

\bigskip
\noindent {\bf Theorem 5. } {\it
Assume that $\sigma \left( \mathbf{x} \right)=\left\{ {{\vartheta }_{1}},{{\vartheta }_{2}},\cdots ,{{\vartheta }_{k}} \right\}$ is of $k$ distinct points. Then there exist vectors $v_1,v_2,\cdots,v_k\in{\mathbb X}$ such that  $x_n=v_1\vartheta_1^n+v_2\vartheta_2^n+\cdots+v_k\vartheta_k^n+o(1)$  as $n\to\infty$. Specially, if $\sigma \left( \mathbf{x} \right)
\subseteq \left\{ 1  \right\}$
 then there exists  $\lim x_n$ as $n\to\infty$.
}

\bigskip
\noindent {\sl Proof: } Let
\[g\left( \lambda  \right)=R\left( \lambda ,\overline{S} \right)\mathbf{\bar{x}}=\sum\limits_{j=1}^{k}{\frac{{{{\mathbf{\bar{c}}}}_{j}}}{\lambda -{{\vartheta }_{j}}}}.\]
This formula holds because the (resolvent)  function $R\left( \lambda ,\overline{S} \right)\mathbf{\bar{x}}$ has simple poles at isolated points $\left\{ {{\vartheta }_{1}},{{\vartheta }_{2}},\cdots ,{{\vartheta }_{k}} \right\}$.
Then
\[   \mathbf{\bar{x}}=\sum\limits_{j=1}^{k}{\frac{\left( \lambda -\overline{S} \right){{{\mathbf{\bar{c}}}}_{j}}}{\lambda -{{\vartheta }_{j}}}}\quad\text{   for every } \lambda \in \mathbb{C}\backslash \left\{ {{\vartheta }_{1}},{{\vartheta }_{2}},\cdots ,{{\vartheta }_{k}} \right\}.
\eqno(*)
\]
Let $\lambda \to \infty $ we get 
$\mathbf{\bar{x}}=\sum\limits_{j=1}^{k}{{{{\mathbf{\bar{c}}}}_{j}}}.$
Replace this back to $\left( * \right)$ we get 
\[\sum\limits_{j=1}^{k}{\frac{\overline{S}{{{\mathbf{\bar{c}}}}_{j}}}{\lambda -{{\vartheta }_{j}}}}=\sum\limits_{j=1}^{k}{\frac{{{\vartheta }_{j}}{{{\mathbf{\bar{c}}}}_{j}}}{\lambda -{{\vartheta }_{j}}}}
\quad\text{  for every  }
\lambda \in \mathbb{C}\backslash \left\{ {{\vartheta }_{1}},{{\vartheta }_{2}},\cdots ,{{\vartheta }_{k}} \right\}
\]
 and consequently, $\overline{S}{{\mathbf{\bar{c}}}_{j}}={{\vartheta }_{j}}{{\mathbf{\bar{c}}}_{j}}$ for $j=1,2,\cdots ,k.$ In the other words,  $\mathbf{\bar{x}}$ is the sum of $k$ eigen-sequences  of the shift operator with respect to $k$ eigenvalues 
${{\vartheta }_{1}},{{\vartheta }_{2}},\cdots ,{{\vartheta }_{k}}$. More exactly, we have
$x_n=v_1\vartheta_1^n+v_2\vartheta_2^n+\cdots+v_k\vartheta_k^n+o(1)$ as $n\to\infty$ where $v_1,v_2,\cdots,v_k\in{\mathbb X}$ are fixed. The proof is now complete.

\bigskip
\bigskip

\sect{ 	Inhomogeneous delay linear difference equations}
\bigskip
\par\noindent 
Now let $B:\mathbb{X}\to \mathbb{X}$ denote  a bounded linear operator. Then $B$ can be extended to the space  ${{\ell }^{\infty }}\left( \mathbb{X} \right)$  by letting ${{\left( B{\mathbf x } \right)}_{n}}=B{{x}_{n}}$ and also to the space $\mathbb{Y}={{\ell }^{\infty }}\left( \mathbb{X} \right)/{{c}_{0}}\left( \mathbb{X} \right)$.
The spectrums of $B$ in the spaces $\mathbb{X}$ and $\mathbb{Y}$ are the same. 
We are interested in the bounded solutions of  the linear difference equation \[{{x}_{n+1}}=B{{x}_{n}}\text{ }+\text{ }{{y}_{n}}\] where $\mathbf{y}=\left\{ {{y}_{0}},{{y}_{1}},\cdots  \right\}$ is a vanishing sequence in ${\mathbb X}$. 
Clearly, for any solution $\mathbf{x}=\left\{ {{x}_{0}},{{x}_{1}},\cdots  \right\}$ we have $S\bar{\mathbf x}=B\bar{\mathbf x}$.
Therefore, the spectrum of any solution $\mathbf{x}=\left\{ {{x}_{0}},{{x}_{1}},\cdots  \right\}$ is contained in the spectrum of the operator $B$ (and in the unit circle).  Theorem 5 gives the following theorem which was proved in \cite{Minh} for the case $k=1$ and $\vartheta_1 =1$.

\bigskip
\noindent {\bf Theorem 6. } {\it  Let $B:\mathbb{X}\to \mathbb{X}$ denote  a bounded linear operator and $\partial \mathbb{D}$ the unit circle. If  $\partial \mathbb{D}\cap \sigma \left( B \right)= \left\{ \vartheta_1,\vartheta_2,\cdots,\vartheta_k  \right\}$ then for every bounded solution $\mathbf{x}=\left\{ {{x}_{0}},{{x}_{1}},\cdots  \right\}$ of
the linear difference equation 
\[{{x}_{n+1}}=B{{x}_{n}}\text{ }+\text{ }{{y}_{n}}
\quad \text{ for } n=0,1,\cdots,
\] 
where $\mathbf{y}=\left\{ {{y}_{0}},{{y}_{1}},\cdots  \right\}$ is a vanishing sequence in ${\mathbb X}$, we have 
$x_n=v_1\vartheta_1^n+v_2\vartheta_2^n+\cdots+v_k\vartheta_k^n+o(1)$  as $n\to\infty$ where $v_1,v_2,\cdots,v_k\in{\mathbb X}$ are fixed. Specially, if $\partial \mathbb{D}\cap \sigma \left( B \right)
\subseteq \left\{ 1  \right\}$
 then there exists  $\lim x_n$ as $n\to\infty$.
}

\bigskip
\noindent For the delay equation \[{{x}_{n+p}}=B{{x}_{n}}\text{ }+\text{ }{{y}_{n}}
\quad \text{ for } n=0,1,\cdots,
\] 
 we have the following result.

\bigskip
\noindent {\bf Theorem 7. } {\it  Let $B:\mathbb{X}\to \mathbb{X}$ denote  a bounded linear operator and $\partial \mathbb{D}$ the unit circle. If  $\partial \mathbb{D}\cap \sigma \left( B \right)\subseteq \left\{ \vartheta  \right\}$ then for every bounded solution $\mathbf{x}=\left\{ {{x}_{0}},{{x}_{1}},\cdots  \right\}$ of
the delay linear  difference equation 
\[{{x}_{n+p}}=B{{x}_{n}}\text{ }+\text{ }{{y}_{n}}
\quad \text{ for } n=0,1,\cdots,
\] 
where $\mathbf{y}=\left\{ {{y}_{0}},{{y}_{1}},\cdots  \right\}$ is a vanishing sequence in ${\mathbb X}$, we have 
 \[\underset{n\to \infty }{\mathop{\lim }}\left({{x}_{n+1}}-\vartheta {{x}_{n}}\right)=0.\]
(Here $p$ denotes a fixed positive integer.)
}

\bigskip

\par\noindent {\sl Proof: }  Clearly, for any bounded solution  $\mathbf{x}=\left\{ {{x}_{0}},{{x}_{1}},\cdots  \right\}$ we have $S^p\bar{\mathbf x}=B\bar{\mathbf x}$.
Therefore, the spectrum of  $S^{p-1}\mathbf{x}$ is contained in the spectrum of the operator $B$ (and in the unit circle).  Consequently, this spectrum is empty or of only one point and our Theorem follows.

\bigskip

\par\noindent {\bf Acknowledgement.}
 Deepest appreciation is extended towards the NAFOSTED  (the National Foundation for Science and Techology Development in Vietnam) for the financial support.


\bigskip

\bigskip


\begin{thebibliography}{99}

\bibitem {Minh} Nguyen Van Minh,  Asymptotic behavior of individual orbits of discrete systems. Proc. Amer. Math. Soc. 137 (2009)  3025-3035. 
\bibitem {Katznelson}  Katznelson Y. and   Tzafriri L., On power bounded operators. J. Funct. Anal. 68 (1986) 313-328.
\bibitem{Rudin} Rudin W. "Real and Complex Analysis"  MacGraw-Hill, New York, 1987.
\bibitem{Vu} Vu Quoc Phong,  A short proof of Y. Katznelson's and L. Tzafriri's theorem. Proc. Amer. Math. Soc. 115 (1992)  1023-1024. 
\end{thebibliography}
\end{document}